\documentclass[12pt]{article}
\usepackage{amssymb, amsfonts, amsmath, amsthm}

\newtheorem{theorem}{Theorem}

\newtheorem{lemma}{Lemma}
\newtheorem{cor}{Corollary}

\newtheorem{remark}{Remark}

\newtheorem{rem}{Remark}

\DeclareMathOperator{\diam}{diam}  
\DeclareMathOperator{\sys}{sys}  
\DeclareMathOperator{\doubl}{doubl}  
\DeclareMathOperator{\length}{length}  

\newcommand{\tr}{\triangle}
\newcommand{\an}{\angle}

\newcommand{\R}{\mathbb R}   
\newcommand{\ghto}{\mathrel{\hbox to 0pt {${}_{_{\,GH}}$\hss}{\longrightarrow}}}

\begin{document}
\title{Bi-Lipschitz equivalent Alexandrov surfaces, II}
\author{Yu.~Burago
\footnote{The  author was partly supported by grants  RFBR
02-01-00090, SS-1914.2003,  CRDF RM1-2381-ST-02, and Shapiro
Foundation of Pennsylvania State Univ. (US).} }
\date{}
\maketitle

\section{Basic definitions and statements}

This paper is a continuation of the paper \cite{BeBu}. Recall that a
map  $f:X \rightarrow Y$ of a metric space  $(X, d_X)$ in a metric
space  $(Y, d_Y)$ is called bi-Lipschitz with a constant  $L$
(or $L$-bi-Lipschitz) if for every $x,y \in X$
$$
L^{-1}d_X(x,y)\le d_Y(f(x),f(y)) \le Ld_X(x,y).
$$
In this case, the spaces  $X$, $Y$ is called  bi-Lipschitz
equivalent (with
constant  $L$). In other words, two metric spaces are
$L$-bi-Lipschitz equivalent  if the Lipschitz distance $d_{Lip}(X,
Y)$ is not greater, than $\lg L$.

Our readers supposed to be familiar with the basic notions of two
dimensional manifolds of bounded total (integral) curvature
theory. Its expositions can be found, for instance,  in
\cite{AlexZalg} and \cite{Reshetnyak}.

 Hereafter the notion of Alexandrov
surface means a {\em complete} two dimensional manifold of bounded
curvature with a boundary; the boundary (which may be empty) is
supposed to consist of a finite number of curves with finite
variation of turn.

{\em Notations:} let  $M$ be an Alexandrov surface with  metric
$d$,\; $\omega$ be its curvature, which is a signed measure,
 $\omega^+$, $\,\omega^-$ be positive and negative parts of the curvature, and
  $\Omega=\omega^++\omega^-$ be  variation of the curvature.
 For any Riemannian manifold $M$ and a Borel set $E\in M$,
$\omega^+(E)=\int_E K^+dS$,\; $\omega^-(E)=\int_E K^-dS$, where $K$ is Gaussian curvature.

A point $p$ carrying curvature $2\pi$ and a boundary point
carrying  turn $\pi$ are called peak points.

We use notation $|xy|$, for the distance  $d(x,y)$;  by
$s(\gamma)$ and $S(E)$ denote the length of a curve  $\gamma$ and the area
of a set $E$, correspondingly.\,  ${\bf D}(X, r)$ means the disk of radius $r$
centered at $X$.

For positive numbers $D,\,C,\,l,\,\epsilon$  and integer $\chi$, by
$\mathfrak M=$\linebreak $\mathfrak M \,(\chi, D, C, l,\epsilon )$ we denote  the class
of closed oriented Alexandrov surfaces $M$ having Euler number $\chi$ and
satisfying the following conditions:

    (i)\; $\diam  M\leq D$,

     (ii)\; $\omega^-(M)\leq C$,

 (iii)\;  if the length of a simple closed curve is less than  $l$, then the curve
 is the boundary of a disk ${\bf D}\subset M$  such that
 $\omega^+({\bf D})\leq 2\pi -\epsilon$.

It follows from (iii) that for every point $p\in M$ the condition
$\omega(p)\leq 2\pi -\epsilon$ holds. Particularly,
$M$ has no  peak points. Besides, the systolic constant for $M$ is not
less than  $l$. (Recall that the systolic constant $\sys M$ of a close surface
$M$ is the infimum of lengths of noncontractible curves in
$M$.

Classes $\mathfrak M$  are compact; the proof is standard, see
Section \ref{section_compact}.

The following theorem is the main result of the paper:
\begin{theorem}
\label{theorem-1}
There exists a positive constant  $L$, depending on
 $\chi, D, C,  l, \epsilon$ only  such that $d_{Lip}(M_1,M_2)\leq L$ for any two
 Alexandrov surfaces  $M_1,\,M_2\in\mathfrak M$.
\end{theorem}

\begin{rem}
\label{rem1}$ $

{\rm
1. A similar theorem is also valid  for nonoriented surfaces.

2. Theorem \ref{theorem-1} is a generalization of Theorem 1 from \cite{BeBu},
but its proof is not independent of the latter one.

3. A generalization of  Theorem \ref{theorem-1} for Alexandrov surfaces with
nonempty boundaries takes place. Naturally, we have to add boundary conditions
in the description of classes $\mathfrak M'$ of surfaces with boundaries. Namely,
distances between boundary components have to be uniformly separated from zero, say
by a number $l$. Besides, for every two boundary points $x,\,y$, the ratio of smaller
boundary arc between $x$ and $y$ to the distance $d(x,y$ also has to be
uniformly separated from zero. The latter condition implies that there is no
boundary point with turn greater than  $\tau (q)=\pi-\epsilon'$ for some fixed
$\epsilon'>0$  (but we do not exclude points $q$ with
 $\tau (q)=-\pi$). We always suppose boundaries to consist of a finite number of
 curves having bounded variation of turn (this condition can be weakened).
We drop the precise formulation because it is a bit
complicated.

In case of surfaces with boundary, the proof is basically the same
as for closed surfaces. Also, it is sufficient to apply Theorem
\ref{theorem-1} to the doubling of a surface with boundary because  freedom in the
choice of a bi-Lipschitz map allows to find it such that it moves boundaries one
to the other.
}
\end{rem}
\medskip

Let $T$ be an end; i.e., an Alexandrov surface
homeomorphic to a closed disk with its center removed and such that
$d(a, p_i)\to\infty$ as $i\to\infty$ for any sequence of points
$p_i\in T$ whose images  in the disk
converge to its center. Here $a$ is a fixed point.
We call the quantity
$v=-\tau(\gamma)-\omega(T)$  the growth speed of the end  $T$. Here
$\tau(\gamma)$ is the turn of the boundary $\gamma$ of the end.
From the Cohn-Vossen inequality it follows that  $v\geq 0$.
Note, that the growth
speed of an end is positive if and only if the limit

$$
v(T)=\lim_{i\to\infty}\frac{l(\gamma_i)}{d(a, p_i)},
$$

\noindent where $l(\gamma_i)$ is the length of the shortest noncontractible loop
with the  vertex $p_i$.
Under condition  $\Omega (T)<\infty$, this limit is well-defined
and is not greater than 2.

Every open (i.e., complete and equipped with an unbounded metric) finitely connected
Alexandrov surface  can be cut (for instance, by geodesic loops) onto
a compact part $M_c$ and ends $T_i$.
Let us consider classes
$\mathfrak  M^* =\mathfrak  M^* (g, C,  l,\epsilon, v_0,  )$
consisting of homeomorphic one to another Alexandrov surfaces $M$ of genus  $g$,
satisfying the conditions
(ii) - (iii) from the definition of class $\mathfrak  M$ and such that all
ends have growth speeds not less, than the number $v_0$
(growth speed of an end does not depend on choice of a loop $\gamma_i$ in its homotopy class).
We will choose loops  $\gamma_i$ in such a way that ends
$T_i$ would satisfy the conditions:
$\Omega(T_i)+\tau^+(\gamma_i)<0,001$, where $\tau$ is the turn from the end side,
and $\sys(\doubl M_c)\geq\sys M_c$. Here $\doubl M_c$ is the double of  $M_c$.
These conditions can definitely be  satisfied if we choose loops
far enough from some fixed point.

Let us denote
$ \tilde{\mathfrak M} =
\tilde{\mathfrak M} (g,D, C, l,\epsilon, s, v_0)$
the subset of class $\mathfrak  M^*$ consisting of  surfaces which can be decomposed
onto a compact part $M_c$ and ends $T_i$ such that the conditions listed above hold true
and, in addition,
  $$
  \diam M_c\leq D,  \quad  \length (\gamma_i)\leq s.
  $$
It is clear that every surface of class
 $\mathfrak  M^*$ belongs to some class
$ \tilde{\mathfrak M}$. Now Corollary below   follows immediately from Theorem
\ref{theorem-1} and Remark 2 from the paper \cite{BeBu}.

\begin{cor} There exists a constant $L_1$, depending on
$g,\,D,\, C, \, l,\, \epsilon,\, s,\, v_0$ only  such that
all Alexandrov surfaces of class
$ \tilde{\mathfrak M} (g,D,  C,  l,\epsilon, s, v_0)$ are
$L_1$-bi-Lipschitz equivalent.
\end{cor}

The author thanks  A.~Belenkiy and V.~Zalgaller whose advises help to
simplify some proofs.
\medskip

{\em A sketch of the proof of Theorem \ref{theorem-1}}

By $L(M,N)$ we denote infimum of Lipschitz constants for
bi-Lipschitz maps $M\to N$, where $M, N\in\mathfrak M (D, C, \chi, l,\epsilon )$.
Suppose that the theorem is not true. Then there exists a sequence of surfaces
$M_i\in\mathfrak M (D, C, \chi, l,\epsilon )$ such that
$L(M_i,N)\to \infty$, where $N$ is a smooth surface of the same class.
It will be shown later that we can suppose surfaces $M_i$  to be equipped with
polyhedral metrics.

Lemma  \ref{lemma_compact} implies that there is a subsequence of  $\{M_i\}$
converging in Gromov--Hausdorff topology and the limit space $M$ for this
subsequence  is an Alexandrov space of the same class
$\mathfrak M (\chi, D, C, l,\epsilon )$. In particular
$\omega (p)\leq 2\pi-\epsilon$ for every point
$p\in M$. Let us   keep  the same notation
for this subsequence.
 From this and Theorem 1 from
\cite{BeBu}, it follows that
$L(M,N)<\infty$. Therefore we come to a contradiction if prove the following lemma.

\begin{lemma}[Key Lemma]
\label{Key-lemma}
Under the assumptions we made above,\\
$L(M_i,M)\leq A<\infty$, where the constant  $A$ does not depend on $i$.
\end{lemma}

The proof of this lemma is the main part of the proof of Theorem 1. It is
exposed in  Section  \ref{sec_key_lemma}.
The proof is based on special triangulations of the surfaces  $M$ and $M_i$ from
Section  \ref{sec_triang} and
on the basic construction of the paper \cite{BonkLang}. Auxiliary
 statements on triangles in $\R^2$ and  Alexandrov surfaces are located
 in Section \ref{sec_triangles}.

\section{Space   $\mathfrak M$ is compact}
\label{section_compact}
\begin{lemma}
\label{lemma_compact}
The space $\mathfrak M =\mathfrak M (\chi, D, C, l,\epsilon )$
is compact in Gromov--Hausdorff topology.
\end{lemma}

{\it Proof.} Precompactness of $\mathfrak M$ was proved in \cite{Sh};
we give here a short proof to make our exposition complete.
Recall, that  $C^*$ means different  constants
depending on parameters of the
class  $\mathfrak M$

1. It is proved in  \cite{Sh} that the space $\mathfrak M$ is precompact.
Nevertheless we give a short proof here to do our text more self-contained.
It is sufficient to show that for any (small enough) $r>0$, on every surface
$M\subset \mathfrak M $, there is a
$r$-net containing not greater than $C^*r^{-2}$ points.

Let us fix   $r<\frac14 l$ and
consider a maximal $2r$-separated set $\{a_1,\dots , a_k\}$ of
points of the surface $M$.
These points form a $4r$-net.
Denote ${\bf D}={\bf D}(a,r)$, where $a=a_i$ and let
$r_0$ be the supremum of numbers $\rho\leq r$ such that the disk
${\bf D}(a, \rho')$ is simply connected for all $\rho'\leq\rho$.

If $r_0\geq\frac12 r\sin\frac{\epsilon}{2}$, then
$S({\bf D})\geq\frac{1}{8}\epsilon r^2(\sin\frac{\epsilon}{2})^2$.
As the whole area of $M$ is not greater than
$(2\pi+\omega^-(M))\diam^2 M\leq D^2(2\pi+C)$, the number of such disks is not bigger than
$C^*r^{-2}$.

Now suppose that  $r_0\leq\frac12 r\sin\frac{\epsilon}{2}$.
Then there is a geodesic loop  ${\bf D}(a, r_0)$
$\gamma$ of length $2r_0$ centered at $a$ separating two components of the
boundary of the disk
${\bf D}(a, r_0)$. As $2r<l$, at least one of components of
$M\setminus \gamma$ being  simply connected. Denote By $K$ its closure.
The Gauss--Bonnet formula says that  $\omega^+(K)\geq\pi$.
The well known inequality for length of a curve in a simply connected region
(see, for example, \cite{Reshetnyak}, section 8.5)
gives
$$
R(K)\leq\frac{2r_0}{\sin\frac{\omega^+(K)}{2}}\leq\frac{2r_0}{\sin\frac{\epsilon}{2}}\leq r,
$$
\noindent where $R(K)$ is inradius of $K$; i.e.,

\centerline{$R(K)=\sup\{ d(x,\partial K), \, x\in K\}$.}

\noindent
This means that $K$  does not intersect disks ${\bf D}(a_j,r)$, $j\not= i$.
Besides, $\omega^+(K)\geq \pi$.
If we add the set $K$ to the disk  ${\bf D}$, then we will get the set which
does not intersect other disks and has curvature
$\geq\pi$. After we perform the same for every disk with radius satisfying the condition
 $r_0\geq\frac12 r\sin\frac{\epsilon}{2}$,
 we get a family of disjoint sets containing our disks.
All different from disks sets  have positive curvature at least $\pi$ each.
Therefore the
number of such sets and the number of all disks can be estimated above
by $C^*r^{-2}$.
\medskip

2. It remains to prove that $M\in\mathfrak M$ if  $M_i\in\mathfrak M$ and $M_i\to  M$.
In  \cite{Sh} it is proved  that $M$ looks like a graph (may be infinite)
some vertices of which ``are blown up''  to Alexandrov surfaces; these surfaces can be
glued together only along separate points, see details in \cite{Sh}.
Therefore it is sufficient to prove that
every point $p\in M$ can not separate its neighborhood $U$.
It becomes clear that  $M\in\mathfrak M$ in this case .
Indeed, obviously $\diam M\leq D$. Curvatures $\omega_i$ of surfaces
 $M_i$  converge weakly (in the sense of K.~Fukaya's definition, see
 \cite{Sh})  to curvature $\omega$  of $M$, therefore $\omega^-(M)\leq C$.
Now it is easy to check that the condition (iii) from the definition of class
$\mathfrak M$ holds for $M$.

So let us prove that any point
$p\in M$ can not separate its neighborhood.
Reasoning to the  contrary,
suppose that there is a point $p\in M$ separating a some its neighborhood. Then
it separates every its smaller neighborhood. Let $p$ separate its round neighborhoods
$U={\bf D}(p, 10r)\supset {\bf D}(p, \rho)={\bf D}$.  Take points
$a$, $b$ in different components of  $U\setminus p$, both at a distance $r$ from $p$.
Let us choose points $p_i,\, a_i,\, b_i\in M_i$ such that
$p_i\ghto p,\;a_i\ghto a,\; b_i\ghto b$ (we mean convergence in the sense of the
Gromov--Hausdorff metric). For all sufficiently big $i$,
 distances $|a_ip|$,  $|b_ip|$ are almost equal to $r$.

Now consider disks $U_i={\bf D}(p_i,10r)$, \,
${\bf D}_i={\bf D}(p_i,\rho)$, where  $\rho\ll r$, for instance
$\rho<\frac{1}{100}(2\pi+C)^{-1}r$ and besides  $r<\frac13 l$.
Note that length of the   disk
${\bf D}_i={\bf D}(p_i, \rho)$ boundary is not greater than $(2\pi+C)\rho$.
As  $r<\frac13 l$, each closed disk  $\bar U_i$,
$\bf{\bar D}_i$ is homeomorphic to an Euclidean closed disk with  not
more than countable set of disjoint open disks removed.

Two cases are possible.

a)  For some subsequence of indexes $i$, the points $a_i$ and $b_i$ are located in
one component of  $U_i\setminus {\bf \bar D}_i$.
 In this case points
$a_i, b_i$  can be connected by a path
of the  length  not greater than
$3r+(2\pi+C)\rho)<4r$ in  $U_i\setminus {\bf \bar D}_i$.
Replace this path by a dotted line with steps
$\frac{1}{10}\rho$ having not bigger, than $40r\rho^{-1}$ points.
Taking the limit, we get a dotted line whose steps are also small  and which ``connects''
 $a$ and $b$ in  $M$.
At least one of the points of this dotted line has to be not farther
than
$\frac{1}{10}\rho$ from $p$.
This contradict to the fact that all distances between points of converging
dotted lines and corresponding points  $p_i$ are not greater than $\rho$.

b) Let points  $a_i$ and $b_i$ be in different components of the set
 $U_i\setminus {\bf \bar D}_i$ (for some subsequence).
In particular, the closed disks
${\bf \bar D}_i$ are not simply connected. Then there is a simple closed loop
in  ${\bf \bar D}_i$ of length not greater than $3\rho$ such that it separates
components containing the points $a_i$ и $b_i$.
This loop is contractible as
$3\rho<l$. Therefore the loop bounds a disk ${\bf D'}$ containing one of our components.
Assume that just $a_i$ are in this component.
The Gauss--Bonnet theorem implies that $\omega^+({\bf D'}\geq\pi$.
Let us choose
$\rho<\frac{1}{100}\sin\frac{\epsilon}{2}(2\pi+C)^{-1}r$.
Then the distance from  $a_i$ to the boundary of  ${\bf D'}$
is not greater than
 $$
\frac{\mbox{boundary length of ${\bf D'}$}}{\sin\frac{\epsilon}{2}}
\leq 3\rho(2\pi+C)(\sin\frac{\epsilon}{2})^{-1}\leq \frac{3}{100}r.
$$
Hence, distances between points $a_i$ and disks ${\bf D}_i$
are not greater than $\frac{3}{100}r$.
Thus, $|pa_i|\leq \rho+\frac{3}{100}r< \frac12 r$.
Contradiction.

The lemma is proved.
\medskip

\section{Lemmas about triangles}
\label{sec_triangles}
Here we collect some auxiliary statements on triangles in Alexandrov surfaces.
These lemmas will be used in Sections 5 and 6. Basically these lemmas are
modifications of  statements proved in  \cite{BeBu} and \cite{AlexZalg}.

Along with usual triangles sometimes we will consider  generalized triangles.
By a generalized triangle, we mean  a disk bounded by three broken lines
(sides of the triangle)
constructed from minimizers. It is supposed that lengths of these sides satisfy
the strict triangle inequality.
We  call total curvature and denote  by $\tilde\Omega(T)$ the
sum of absolute curvature of a
generalized triangle $T=\tr ABC$ and variations of turn of its sides; i.e.,
$\tilde\Omega(T)=\Omega(T)+\sigma (AB)+\sigma (BC)+\sigma (CA)$,
where $\sigma$ means variation of turn from the triangle side.
The angles of a generalized triangle are allowed to be zero.
For short we will  drop sometimes the word ``generalized''.

Recall, that a comparison triangle for a (generalized) triangle $T$ in Alexandrov space
$M$ is a planar triangle with the same  side lengths.

Usually we will consider generalized triangles $T$ for which
$\tilde\Omega(T)$ is small enough. If this quantity is small in comparision with
the angles of a triangle, then such a generalized triangle is bi-Lipschitz
equivalent to its comparison triangle, where Lipschitz constant depends on
low angles estimate. More precisely, the following statement takes place.

\begin{lemma}
\label{small-tr}
For any $\alpha>0$,  $L>1$, there exists  $\delta=\delta(\alpha, L)>0$
with the following property.
If every angle of a generalized triangle $\tr ABC$ is not less than
$\alpha$ and  $\tilde\Omega(\tr ABC)<\delta$, then there exists a $L$-bi-Lipschitz
map of the  generalized triangle $\tr ABC$ onto its comparison triangle,
this map may be chosen in such a way that its restriction on the boundary of the
triangle is an isometry which moves every vertex to a vertex.
\end{lemma}

This lemma is a minor modification of Lemma 4 from
 \cite{BeBu} and can be proved by the  same way. By this reason we drop the proof.
We will also need a more general statement.

\begin{lemma}
\label{small_tr2}
Let a simply connected closed region $T$ is equipped with a polyhedral metric and
bounded by two shortest curves $BA,\;BC$ and a geodesic broken line $AC$.
Suppose that $|BA|+|BC|>s(AC)$, where $s(AC)$ is the length of $AC$.

Assume that  $T$ is starlike with respect to a point $C$; i.e., all shortest curves
$BX$, where $X\in AC$, intersect $AC$ at point
$X$ only. Let angle   $\an ABC$ satisfy the condition $0<\phi\leq\an ABC\leq\frac{1}{10}$,
Also suppose that for every $X\in AC$,  angles between a shortest $BX$
and started at $X$ arcs of
the broken line $AC$ are in the interval
$[\frac{\pi}{2}-\frac{1}{10}, \frac{\pi}{2}+\frac{1}{10}]$.

Then there exist constants $\delta$, $L$
such that if   $\tilde\Omega(T)<\delta$,  then  $T$ is  $L$-bi-Lipschitz equivalent
to a planar triangle
$\tr A'B'C'$, whose side lengths are equal to
$|AB|,\,|CB|,\,s(AC)$, correspondingly.

If in addition  $AC$ is a shortest curve, then one can choose
$L$ as a function
$L=L(\delta)$ in such a way that $L\to 1$ as $\delta\to 0$.
\end{lemma}

Recall that we suppose that $L$-bi-Lipschitz map of
$T$ onto its  ``comparison triangle''  $A'B'C'$
keeps lengths of boundary curves fixed.

This lemma also is a modification of Lemma 4 from \cite{BeBu},
and can be proved by the same way, so we omit details of the proof.
The idea of the proof is the following. First of all
we map  $T$ onto a planar closed region $\tilde T$ bounded by intervals $A_1B_1$, $C_1B_1$
and a broken line  $A_1C_1$  such that  $|A_1B_1|=|AB|$,\,
$|C_1B_1|=|CB|$,\, $\an A_1B_1C_1=\an ABC$. To do this we use  Tchebyshev coordinates.
One can verify  that the turns of the broken line
 $A_1C_1$  at its vertices can be estimated above by some value depending on
 $\phi,\;\phi_1$ and smallness of $\delta$. After that, it is not difficult to map
  $\tilde T$ onto the comparison triangle   $\tr A'B'C'$.

Besides we will need the  following corollary.
\begin{cor}
\label{corr-4-angle}
Let a quadrangle $\Box=AA_1C_1C$ be boundary convex  and bounded by
four shortest curves,  $\tilde\Omega(\Box)<\delta$.
Suppose that

\centerline{$|AA_1|=|CC_1|$,\; $|AC|<\frac12|AA_1|$,\; $|A_1C_1|<\frac14|AA_1|,$}

\centerline{
 $\an A_1<\frac{\pi}{2}-\phi$,\;  $\an C_1<\frac{\pi}{2}-\phi$,\;
 $|\an A-\frac{\pi}{2}|<\phi$, \; $|\an C-\frac{\pi}{2}|<\phi$,}

\noindent where $0<\phi<\frac{1}{10}$.
Then for every fixed $\phi$, there is a function
$L=L(\delta)\geq 1$  such that
$L\to 1$ and $g\to 0$ as  $\delta\to 0$ and
$\Box$ is $L$-bi-Lipschitz equivalent to a planar  quadrangle having the same
side lengths and satisfying the condition: differences between its angles
$\an A'_1$, $\an C'_1$ and angles $\an A_1$, $\an C_1$ are not greater than
$C^*\delta$.
\end{cor}

To prove let us separate the  quadrangle $\Box$ from the surface
and attach a planar triangle $ОA_1C_1$ along
$A_1C_1$  such that its sides $A_1O$, $C_1O$
are continuations of the  quadrangle sides; i.e.,  they form angles $\pi$
with the shortest curves $A_1A$,   $C_1A_C$, correspondingly. Thus, we obtain a
generalized triangle $T=\tr OAC$ (it is not necessary an ordinary triangle as its
sides can be not shortest curves). It is not difficult to check that this triangle
satisfies
the conditions of Lemma  \ref{small_tr2}. Applying  this lemma gives a
bi-Lipschitz (with a constant depending on $\phi$ and
smallness of  $\delta$ only) map
$f_0\colon G\to \tr O'A'C'$, where
$\tr O'A'C'$  is a comparison triangle for $T$;
restrictions of $f_0$ on the sides are isometries.

In the proof of Lemma  \ref{small_tr2}, the map $f_0$  is constructed in two steps.
First we map $T$ onto a planar figure
bounded by two intervals (the images of $OA$ and $OC$) and the broken line
$\gamma$ (the image of the shortest curve $AC$).
To do this we use Tchebyshev coordinates. As the second step, the broken line
$\gamma$ is transformed into an interval, see details in \cite{BeBu}.

 As triangle   $\tr OA_1C_1$ is planar, the first map acts isometrically on it,
in particular, the shortest curve  $A_1C_1$ is mapped onto
 interval $A'_1C'_1$ of the same length. Now it is not difficult to
 straighten up
  the broken line  $\gamma$ keeping interval  $A'_1C'_1$ fixed. To do this
 let us cut the quadrangle   $A'A'_1C'_1C'$ by the diagonal $A'_1C'$
 into two triangles. Now we can straighten  the broken line   $\gamma$ as a side of
 the ``curved triangle'' $\tr A'_1C'A'$. For this we transform  $\tr A'_1C'A'$
 the same way as it was done in the item
 8 of the proof of Lemma  4
 in \cite{BeBu}. We keep the
 triangle $A'_1C'_1C'$ firm during this process.

\begin{remark}
\label{boundary_bilipsch}{\rm
The words ``bi-Lipschitz equivalence'' will always mean (if contrary is not supposed)
the existence of a bi-Lipschitz map with a constant depending on parameters of the
 class  $\mathfrak M$ only. If a surface has the boundary, we suppose that the
restriction of a bi-Lipschitz map on the boundary is linear. In case
of triangles we also suppose that  vertices are mapped into vertices.
}
\end{remark}

The total curvature  $\tilde\Omega (G)$ of a subset  $G$ of a generalized triangle
$T=\tr ABC$
is equal, by definition, to the sum of $\Omega(G)$
and negative  turn of intersection of triangle sides
with $G$ (we mean open sides without vertices). Recall that turn of an ordinary
triangle side is nonpositive.)

By shortest curves connecting points of a triangle we mean shortest curves of its
induced metric.

\begin{lemma}
\label{concentration}
For any positive  $\Psi,\,R,\,\delta$, there exists a number  $r>0$ having the following
properties.
Let a simple triangle  $\tr ABC$ satisfy the conditions:
$\tilde\Omega (\tr ABC)<\Psi$, $\tilde\Omega (\tr ABC\setminus {\bf D}(Z,r))<\delta$.
Then,

(i) if $B=Z$, $|AC|<R$, $d(B, [AC]>R$, then the differences between angles
$\an B,\;\an C$ of the triangle and corresponding
angles of its comparison triangle $\tr A'B'C'$ are not greater than $2\delta$.

(ii) If  $A,\;B\in {\bf D}(Z,r)$, and
$|CZ|\geq R$, then
$\an ACB-\an A'C'B'\leq 2\delta$, where $\an A'C'B'$ is
the angle in the comparison triangle.
\end{lemma}

\begin{remark}
\label{remark3}
 If we choose  $r$ such that
$\an A'C'B'<\delta$ in the item (ii), then obviously
$\an ACB<3\delta$.
\end{remark}

Here we restrict ourselves  by a sketch of a
proof, because technique of the proof
is the same as in section 2 of chapter IV in the book
\cite{AlexZalg}; the reader can find all details in the book.
(Note, that it is enough to prove the lemma for polyhedral metrics only;
by the way,  we need only this case.)

In the item (i), the idea of the proof is the following: suppose that in our triangle
(with a polyhedral metric),  there are
 points of positive curvature  at the distance less than
$r<R/2$ from $B$. Then one can consecutively move these points $X$ until
they are placed at the
distance at least $R/2$ from $B$. For this we look for a bigon
(bounded
by two shortest lines with common ends at $B$ and one more point $Y$) containing
point $X$ and then remove the bigon.
As a result, vertex $X$ vanishes but
additional curvature can appear
at the point $Y$. This additional curvature at least $\frac{2r}{R}$ times less
than curvature of the removed vertex $X$.
This means that curvature of the vertex will be less
than $\delta/2$ if $4\Psi r< \delta R$.

Now there is no positive curvature in  $R/2$-neighborhood of $B$. This allows
to move all vertices $X$ of negative curvature
at the distance at least $R/2$ from $B$.
To do this we glue an additional material in a slit looking like a tree with one vertex;
it consists of  $BX$ and several additional slits started at  $X$.
At this step negative curvature decreases  almost in the same
proportion as positive curvature has been decreased. As a result variation of curvature
becomes less than $2\delta$. The side $AC$ keeps to be a shortest
during this process because it was far enough of the deformed region of the triangle.
Angles  $\an A,\;\an C$ were not changed too. This proves the item (i).
\medskip

In the item (ii) the idea of the proof is almost the same: at the
first step we remove all vertices of positive curvature on the
side $AB$ by cutting  bigons with vertex $C$. This allows to
remove all positive curvature. Choosing $r$ as in the item (i)
we can guarantee that change of angle $\an C$ is not greater than
$\delta$. However the side  $AB$ can cease to be a shortest curve.
Let us replace it in such a case by a shortest curve (in the
induced metric), which is not longer. As a result, variation of
curvature can only decrease. Applying the angle comparison theorem
to the triangle of nonpositive curvature, bounded by $AC,\; BC$
and a new shortest curve $AB$ immediately gives the required
inequality.

\section{Approximations and triangulations}
\label{sec_triang}
\begin{lemma}
\label{approx-2}
Every compact Alexandrov surface $M$ (possibly with boundary) without peak points
can be Lipschitz approximated by surfaces $P_i$ with polyhedral metrics.
Moreover,   convergence $P_i\to M$ can be made  regular; the latter means that
$\omega_i^\pm\stackrel{weak}{\longrightarrow} \omega ^\pm$.
\end{lemma}

This lemma was announced by Yu.~Reshetnyak in \cite{Reshetnyak1}
(actually in a more general form), but the proof has never been published.

{\em Proof.} Recall that a triangle is simple if it is boundary convex,
its sides have no common points except vertices and bound a disk.
According to \cite{AlexZalg}, Theorem 3 of Chapter 3,
$M$ can be partitioned onto arbitrary small simple triangles such that all
triangle inequalities are strict. In addition, for any finite set of points and
a finite set of shortest lines started at these points, it is possible to
include these points
to the set of vertices and some initial intervals of the shortest lines
to the set of edges.
Replacing each triangle of the partition by a planar triangle with the same
side lengths (comparison triangle), we get a surface  $P$ equipped with a polyhedral
metric. It is proved in \cite{AlexZalg}, Theorem 7 of Chapter 7, that if
triangles of the partitions become smaller and smaller, the
sequence of  polyhedra $P_i$ converges to $M$ uniformly and regularly.

Now we particularize our partition according to
the purpose to provide Lipschitz convergence.
Namely, let
$\theta_0=\frac{1}{100}\min_{p\in M}(2\pi-\omega(p)$.
There is only a finite number of points with absolute curvatures
greater than $\theta_0$. Denote them by $E_1,\dots , E_m$.
We construct a partition such that the star of each point $E_k$, $k=1,\dots , m$,
consists of  isosceles triangles with vertex   $E_k$, angles of the triangles at
$E_k$ being in the interval
$(2\theta_0,\, 10\theta_0)$. Besides, we do triangles of the partition so small
that $\tilde\Omega (T)<0,001\theta_0$ for every triangle  $T$.
As curvature of triangles is small, all the angles except
may be one angle in any triangle to be
less than $\pi -5\theta_0$. After we cut each triangle with a ``big'' angle
onto two triangles we get a partition such that all angles of triangles are less
than $\pi -5\theta_0$.

Now we change slightly our partition to get a partition  all angles of which
are positive. To do this we  replace some ordinary triangles by generalized ones.
We can  do this in such a way that the stars of points $E_k$ do not change
and every changed side is transformed to a broken geodesic having almost
the same length and turn as the replaced side (see details in
Lemma 6 of the paper \cite{BeBu}). This deformation is supposed
to be so small that all the properties listed above are preserved.

Let $i$ be so great that $\frac1i\ll\theta_0$.  By
$\delta_i=\delta(\theta_0,\,\frac1i)>0$ denote the number corresponding to
$\theta_0$ and
$L_i=1+\frac1i$ in according with Lemma \ref{small-tr}.
We can choose the partition of $M$ onto generalized triangles
$T$ so that (in addition to properties mentioned above) the following holds:
$\tilde\Omega(T)< 0,001\min\{\delta_i,\,\theta_0\}$ and $\diam T<\frac1i$.

By $M^i$ we denote the surface   $M$ jointly with the partition we have chosen.
All angles of (generalized) triangles $T_{ij}$ of this partition are not
zero and, therefore, they are not less than some number
$\beta_i>0$.
The  triangles $T_{ij}$ having all angles not less than
$\theta_0$ are $L_i$-bi-Lipschitz equivalent to
their comparison triangles (Lemma \ref{small-tr}).
In particular, it takes place for all triangles adjacent to vertices $E_k$.

Only one angle of any triangle  $T_{ij}$  can be less than $\theta_0$
(because the triangle does not have ``big'' angles and its curvature  is small).
Let  angle $\an ABC$ of triangle
$\tr ABC$ be less than $\theta_0$ and its other angles be greater than $\theta_0$.
Obviously such a triangle can not be adjacent to  vertices $E_k$.
Take points $A_1$, $C_1$, $B_1$ on the sides $AB$, $BC$, $AC$ so that
$|AA_1|=|AB_1|,\; |CC_1|=|CB_1|,\; |BA_1|=|BC_1|$ (``Gromov's product'').
Let us connect these points with shortest lines
in the induced metric of the triangle.
Note that due to smallness of $\tilde\Omega(\tr ABC)$,
these shortest lines will cut
$\tr ABC$ onto 4 simple (generalized)triangles, all angles of these triangles,
except may be $\an A_1BC_1$, being greater than $\theta_0$.
Now we choose points
$A_2\in A_1B$, $C_2\in C_1B$, $B_2\in A_1C_1$ such that
$|A_1A_2|=|A_1B_2|,\; |C_1C_2|=|C_1B_2|,\; |BA_2|=|BC_2|$ and continue this process.
It is not  difficult to calculate that, as curvature is small, all
angles of triangles
 $\tr A_kA_{k+1}B_{k+1}$,
 $\tr C_kC_{k+1}B_{k+1}$,  $\tr B_kA_kC_k$, are bounded below by  $\theta_0$
 (we set  $A=A_0,\; C=C_0$, $k=0, 1, \dots$) and for sides of these triangles the
 strict triangle inequality holds. Hence, all these triangles are
$L_i$-bi-Lipschitz equivalent to their comparison triangles.
It is easy to see that $A_k\to B,\; B_k\to B$ as
$k\to\infty$. Therefore there is a number $k$ such that
$\tilde\Omega(\tr A_kBC_k)<\delta(\beta_i, L_i)$. This means that
$\tr A_kBC_k$ \; is  $L_i$-bi-Lipschitz equivalent to its
 comparison triangle (Lemma \ref{small-tr}).
Now, replacing each triangle of our partition of $\tr ABC$ by its comparison triangle
(and, of course, doing this for each triangle  $\tr ABC$) we obtain a polyhedron
$P_i$, which is
$L_i$-bi-Lipschitz equivalent to $M$.

 Lemma \ref{approx-2}  is proved.
 \medskip

\begin{lemma}
\label{triangulation} For every $\nu>0,\; d>0$,
each compact (possibly with boundary) Alexandrov surface $M$
without peak points has a triangulation
$\{T_k\}$ such that

(i) $\tilde\Omega(T_k)<\nu$;

(ii) $\diam T_k<d$;

(iii) all angles of triangles $T_k$ are not less than $\alpha(\theta)>0$
where $\alpha$ depends on
$\theta=\min\{\min\{2\pi-\omega^+(p)\colon\; p\in M\},
\min\{\pi-\tau (q)\colon\, q\in\partial M\}\}$ only. Here
 $\tau(q)$ is turn at point  $q$.

(iv)  The set of vertices contains any given a priori finite set of points $E_k\in M$.
\end{lemma}

\begin{remark}
\label{remark} {\rm
a) In case the metric of  $M$ is polyhedral, this lemma was in fact proved
in \cite{Bur} (Theorem  2) on the basis of the theorem from  \cite{BurZalg}
(see also \cite{BurZalg1}); all triangles of the
triangulation are flat in this special case.

b)  Probably, using Tchebyshev coordinate, it is possible to prove
Lemma \ref{triangulation} the
same way as it has been proved for polyhedra in \cite{Bur}, \cite{BurZalg}.
 However it is simpler to reduce Lemma  \ref{triangulation} to the case
 of polyhedra with the help of Lemma \ref{approx-2}.
}
\end{remark}

{\em Proof.}
From Lemma \ref{approx-2} it follows that
$M$ can be Lipschitz approximated
by polyhedra  $P_i$. Let $f\colon P_i\to M$ be corresponding
$L_i$-bi-Lipschitz maps, $L_i\to 1$ as $i\to\infty$. Fix a  set  $\{F_k\}$ in $M$.
We include  all the points with variation of  curvature
greater than  $\frac{1}{10}\nu$ to this set. Denote $F_{ki}=f_i^{-1}(F_k)$.
As it was mentioned, $P_i$ can be triangulated onto planar triangles satisfying
conditions  (i)-(iv) of the lemma, even if we replace numbers $\nu$, $d$ to
$\frac{1}{100}\nu$, $\frac{1}{10}d$ beforehand.
Choosing such a triangulation of $P_i$ we can include all points $F_{ki}$ to the set
of vertices. Also we can suppose the triangles to be so small that every
$d$-neighborhood of each point $A\in M$ contains not more than one point
$F_k$ and absolute curvature of such a neighborhood without point
$F_k$ is not greater than $\frac{1}{20}\nu$. Also we can suppose that
the similar is true
for every polyhedron  $P_i$ if  $i$ is big enough. One can choose the
described triangulation of the polyhedron $P_i$ in such a way that all angles
of the triangles are bounded below by some number  $2\alpha$ depending on
$\theta_i=\min\{\min\{2\pi-\omega^+(p)\colon\; p\in P_i\},
\min\{\tau (q)\colon\, q\in\partial P_i\}\}$ only; in particular,
$2\alpha$ does not depend on smallness of triangles.
(Note, that numbers $\theta_i$ for polyhedra
$P_i$ with great $i$ are almost the same as the corresponding number $\theta$ for $M$.)
Let us set $\alpha$ to be equal a half of this number. Now we use Lemma  \ref{small-tr}.
As  $\alpha$ does not depend on smallness of triangles, the choice of points $(F_k)$
and numbers $\nu$, $d$, we can assume  $\nu$ to be
so small in comparison with $\alpha$
that $2\nu<\delta=\delta(\alpha, L=2)$, where  $\delta$ is defined
by Lemma \ref{small-tr}.

Now connect by shortest curves points of $M$, whose inverse images in $P_i$
are connected by shortest curves (keeping $i$ fixed).
We claim that, if $i$ is big enough, this makes a triangulation of $M$
combinatorially equivalent to the triangulation of  $P_i$,
all angles of this triangulation being separated from zero by a number
depending on $\theta$ only and the angles at  $F_k$ being only slightly
(less than  $2\nu$) different from corresponding angles at  $F_k$.

Indeed, let $AB$ and $BC$ be the edges of the triangulation of $P_i$,
$A'B'$ и $B'C'$  the shortest curves in $M$, correspondingly.
The shortest curves $AB$ и $BC$   divide a neighborhood of $B$ onto two sectors.
The sector corresponding to triangle $ABC$ is distinctly smaller and its angle
is equal to the angle $\an ABC$ of the triangle. In addition
$\an A'B'C'$ is
almost equal to $\an ABC$ if $i$ is big.
 Combinatorial equivalence of the nets follows easily from this.
Other properties of the triangulation of $M$ now follow from corresponding properties
of triangulations of polyhedra $P_i$ (if $i$ is sufficiently big).

Lemma \ref{triangulation} is proved.

\section{Proof of Key Lemma}
\label{sec_key_lemma}

1. {\em Preliminary agreements.}  Here we will consider only a sequence
of surfaces  $M_j\in\mathfrak M$ converging (in Gromov--Hausdorff topology) to
a surface  $M$.  We will construct partitions of these surfaces into triangles.
These triangles we suppose to be so small that the values of arguments
$\chi,\,D,\,l\,$
of class  $\mathfrak M$ do not play any role in our consideration.
By bi-Lipschitz equivalence of  triangles or more general figures,
we  {\it always} mean a bi-Lipschitz map with a constant depending on
$C$ and $\epsilon$ only.
If there are marked points in the boundary of a figure
(we claim that vertices of a triangle
are always  marked), we assume that our map moves
marked points to marked ones and that the restriction of the map on
boundary curves connecting
marked points is linear.
\medskip

2. {\em Choice of scales.} We have three scales. First, it is the
size of angles of triangles.
Partitions of the limit surface  $M$ are constructed
 of two types triangles: ``ordinary'' and
``special'' ones.
In accordance with Lemma  \ref{triangulation}, angles of ordinary triangles are
separated from zero by some constant  $\lambda>0$ depending on  $C$ and
$\epsilon$ only. All special triangles are isosceles, and angles at their vertices
belong to the interval  $(\varphi_0, \varphi_1)$ where
$\varphi_i$ are small positive numbers also  depending on  $C$ and
$\epsilon$ only; they will be chosen in item 3 of the proof.

At the second step we choose a positive number $\delta$ to be so small that
conclusions of Lemma \ref{small-tr} and Lemma \ref{small_tr2} holds even if
angles of triangles are bounded below by the number
 $0,01\pi\varphi_0(2\pi+C)^{-1}$ instead of $\varphi_0$. Some
 quantities such that they can be estimated
 above by  $C^*\delta$, where $C^*$ depends on $C$ and
$\epsilon$ only, will arise in the process of  the
 proof.  By  $\delta'$ we denote such quantities. It is important that we can
 unboundedly decrease $\delta$ and, therefore,  $\delta'$ keeping
 $C$ and $\epsilon$ fixed. By this reason we will drop a factor $m$ in quantities
 of the form $m\delta$ if $m$ is not too big (say, less than 50).
 It is convenient to assume that
 $\delta'\ll\min\{\varphi_0,\lambda\}$.

 After we have fixed $\varphi_i$ and $\delta$ we  choose a partition of
 $M$ into so small triangles  that variation of  curvature for every triangle is
 less than  $<\delta$. (By variation of  curvature for a triangle $T$ we mean
 $\tilde\Omega(T)$.) In fact we choose
 the partition even more petty. This helps us to
transfer the partition to the surfaces $M_j$ for  big values of $j$.

 Finally, fixing  a partition, we choose so great integer  $j_0$, that
 for $j>j_0$ essential portions of curvature of $M_j$ are concentrated
 in very small (in comparison with size of the triangles) neighborhoods
 of vertices.

 Let us explain the last point. K.~Fukaya defined  weak convergence
 of measures for the case of Gromov--Hausdorff convergence
 of spaces,  see details in \cite{Sh}.
 For a subsequence curvatures $\omega_j$  of $M_j$
 converge weakly to curvature $\omega$ of $M$; positive and negative parts
$\omega_j^+,\;\omega_j^-$ of  $\omega$ converge weakly to some finite measures
$\mu^+,\;\mu^-$. We have
$\mu^\pm\geq\omega^\pm$, where $\omega^\pm$  are positive
and negative parts of  $\omega$. Choosing a partition of $M$ onto triangles
we require that not only variation of curvature  but also  measures
$\mu^+$ и $\mu^-$ be small (less than $\delta$) on all triangles with
 vertices removed. (Note, that  both measures, $\mu^+$ and $\mu^-$,
can be big simultaneously  at a vertex. The reason is that the  convergence $M_j\to M$
can be  nonregular.  All vertices for which  these measures
are big are  special.) However for converging surfaces
 $M_j$, measures $\omega^\pm$ are not necessary concentrated at vertices,
 they can be ``spread out''. Hopefully, for any $R>0$ there exists a number
$j_0$ such that, for any vertex  $B$ of a special triangle $T$, almost all
$\omega^\pm(T)$ are concentrated in $R$-neighborhood $V={\bf D}(B,R)$
of point $B_j$ for $j>j_0$.

 This means that for every special triangle  $\tr A_jB_jC_j=T$
\begin{equation}
\label{ineq-concentration}
\omega_j^\pm (T\setminus V)<\delta.
\end{equation}
later on we suppose   $j$ to be so big that the inequality
\eqref{ineq-concentration}
holds true for  $R$ we have chosen.
\medskip

3. {\em Special vertices and triangles.}
Let us consider the limit surface  $M$. Its partition will be
based on Lemma  \ref{triangulation}.
Before applying the lemma we  choose a finite set of points $F_k\in M$
and triangulate small closed neighborhoods $ Q^0_k$ of these points in a special way.
We suppose that
$ Q^0_k\cap Q^0_l=\emptyset$ for $k\not= l$. After that,
 we apply our Lemma \ref{triangulation} to the surface
$M_0=M\setminus \cup Q^0_k$ with boundary. As a result  we obtain
a partition of $M$ onto triangles.

We set $Q^0_k$ to be stars of points $F_k$. These stars consist of isosceles triangles
$\tr F_kA_{ki}A_{k\,i+1}$, where $|F_kA_{ki}|=|F_kA_{k\,i+1}|$.
We call points $F_k$ and triangles
$\tr F_kA_{ki}A_{k\,i+1}$ adjacent to them to be {\em special}.
Construction of these stars has some freedom; in particular
 angles of the special triangles, their size and pettiness of triangulation
can be changed.
We will use this freedom as follows.

Let $C$ and $\epsilon$ be constants from the definition of class $\mathfrak M$,
$C_1=2\pi +C$.
First we choose intervals for values of the angles with vertices at $F_k$
(before choosing  points  $F_k$). These angles should be so small that even
being multiplied by
 $2\pi/\epsilon$ they remain ``small'', say, less than $0,001$.
From the other hand we should bound uniformly these angles below and bound a
number of edges
at a special vertex above.  So we require that
these angles $\psi$ to be in the interval
\begin{equation}
\label{size-angles}
\varphi_0=10^{-5}{\epsilon}<\psi <\varphi_1=10^{-4}\epsilon.
\end{equation}

These conditions are always  met in such a way that
the number $m$ of edges at $F_k$ is uniformly bounded above:
\begin{equation}
\label{number-edges}
m<10^6(2\pi+C)\epsilon^{-1}.
\end{equation}

Now we choose the
number $\delta$, which characterizes smallness of curvature of triangles.
Namely, set $L=\frac{11}{10}$  and let $\delta_1$ be a number corresponding to
the numbers $L$ and $\alpha=\lambda$ in according with Lemma \ref{small-tr}.
Similarly we can find $\delta_2$, corresponding to
$L$ и $0,01\varphi_0C_1^{-1}$.
Then Lemma  \ref{small_tr2} gives us   $\delta_3$,  corresponding to
$\phi=0,1\varphi_0$, $\phi_1=0,01$.
Finally we put $\delta =\frac{1}{100}\min_i\{\delta_i\}$.
Hence,  $\delta$ depends on $C$, $\epsilon$ only.
Recall that we can decrease $\delta$ if we need and after that find a partition
of $M$ onto triangles such that absolute curvatures of the triangles do not exceed the
new value of $\delta$; low bounds of triangle angles will not be changed.
As $\varphi_1<0,01$, we can suppose that adjacent to the base  angles of
special triangles are close to $\pi/2$ (up to $\varphi_0$).

After we fix  set $\{F_k\}$ (we will do that some later) we will choose
stars  $Q^0_k$ of these vertices to be so small that
$\Omega (Q^0_k\setminus F_k)<\delta$ and besides $\diam Q^0_k<\delta$.
Hence, each special triangle will be
$\frac{11}{10}$-bi-Lipschitz equivalent to its comparison triangle.
Note that turn of the boundary of  $Q^0_k$ from outside at any point is not big,
 say, less than  $\pi/2$.

4. {\em Partition of $M$ onto triangles.}
Let us triangulate the surface  $M'=M\setminus \cup_kQ^0_k$
in according with Lemma  \ref{triangulation}.  All angles  of such a triangulation are
bounded below by some number $\lambda>0$ depending on
the number $\theta$ of $M_0$ (see item (iii) of Lemma \ref{triangulation}).
The last number actually does not depend on our choice of vertices  $F_k$ and their stars,
so we can set $\theta=\epsilon$. Indeed, as it was mentioned above, outside turn of
the boundary of any star at any point is not greater than
 $\frac12\pi$.
At the same time, including  all points having big values of $\mu^\pm$ in the set $\{F_k\}$,
and taking a sufficient petty triangulation, we can provide the inequality
$\tilde\Omega (T)<\delta $
for all triangles, with  $\delta$ as chosen above.

Thus, from the beginning we include all the points having
curvature $\Omega(F_k)\geq\delta$ in  $F_k$; after that we choose
stars  $Q_k$ to be so small  that $\Omega(Q_k\setminus F_k)\leq\delta$;
and finally we triangulate  $M_0$ so that
for any triangle  $Т$ the inequality $\tilde\Omega T\leq\delta$ holds.
This is possible, as our constants do not depend on the choice of the set of points
$F_k$, stars $Q_k$ and a triangulation. In fact, we will add  some requirements
(which can easily be fulfilled) on the choice
of partition of $M$ in the beginning of item 5.

As a result we get a partition of $M$ onto two kinds of triangles: special  ones and
others, each triangle $T$ satisfying
$\tilde\Omega (T)<\delta$ and being
$\frac{11}{10}$-bi-Lipschitz equivalent to its
comparison triangle.
\medskip

5. {\em Converging surfaces.}
Lemma \ref{approx-2} allows us to think that  converging surfaces $M_j$ are
equipped with polyhedral metrics. Taking a subsequence, we can suppose that
curvatures $\omega_j$ of surfaces  $M_j$ converge weakly
(in the sense of definition from  \cite{Sh}) to curvature of  $M$,
their positive and negative parts  $\omega_j^+,\;\omega_j^-$ converge weakly
to some finite measures
$\mu^+,\;\mu^-$. Recall that  $\mu^\pm\geq\omega^\pm$, where $\omega^\pm$ are
positive and negative parts of curvature of  $M$.

Consider the partition of $M$ chosen in the  item 4 of the proof.
Let $\{A_i\}$ be the set of all the vertices of the partition,
$\{F_k\}$ be its subset consisting of the special vertices.
Taking more reach set
$\{F_k\}$, small stars $Q_k$  and making triangles smaller, we can
include  all points
$X\in M$ with
$\mu^\pm (X)\geq \delta$ in set  $\{F_k\}$ and ensure every
closed triangle with vertices removed to satisfy the inequality
$\mu^\pm<\delta$.
The condition (iii) from definition of classes
$\mathfrak M$ implies $\mu^+(X)\leq 2\pi-\epsilon$ for every point
$X\in M$. As a result, we can ensure all triangles to be so small that
the inequality
$\mu^+(Q^0_k)<2\pi-\frac23\epsilon$ holds for each star.

Let diameters of all triangles are not greater than a number $d>0$ so small that

(a)  the inequality $\mu^\pm(E)<\delta$  holds for every set $E$
such that it does not contain points $F_k$ and its diameter
$\diam (E)\leq 10d$;

(b)  each circle of radius $10d$ contains not more than one vertex $F_k$.

Denote by $ A_{jk}$  points of the surface
 $M_j$ such that $A_{jk}\ghto A_k$  as $j\to\infty$;  in particular, $F_{jk}\ghto F_k$.
For a vertex $A_k$ belonging to the boundary of a star $Q^0_{l}$, let us  choose
points  $A_{jk}$ so that
$|F_{jl}A_{jk}|=|F_lA_k|$.

Later on we assume numbers $j$ to be so big that if a set $B\subset M_j$ has
diameter $\leq 6d$ and does not intersect $\delta$-neighborhoods of points $F_{jk}$,
then
$\Omega_j(B)=\omega_j^+(B)+\omega_j^-(B)<\delta$
\medskip

6. {\em Partitions of surfaces  $M_j$ and non-special triangles.}
To construct a partition of the  surfaces $M_j$, connect pairs of
points $A_{jk}$ by shortest curves if and only if corresponding pairs of points
 $A_k$ are connected by shortest curves. Such shortest curves are not
 necessary unique and can have superfluous intersections one with another. We
 will choose shortest curves in  a way to avoid such extra intersections.
 Note, that shortest curves connecting  $A_{jk}$ with $A_{js}$ are not necessary converge
 (in Gromov--Hausdorff  metric sense) to  shortest paths between
  $A_k$ с $A_s$ chosen beforehand.
Almost the same arguments as in Lemma  \ref{triangulation} show that we get a partition
combinatorial equivalent to the partition of the surface $M$.

 Let  $\tr ABC$ of the surface $M$ be non-special. Its angles are almost the same as
 angles of its  comparison triangle. If numbers  $j$ are great enough, triangles
 $\tr A_jB_jC_j$ are in  regions with small variation of curvature
 (less than $\delta$). Hence, the angles of such a triangle are almost equal
 to the angles
 of its comparison triangle.  Lemma 4 from  \cite{BeBu} implies that  both triangles,
$\tr ABC$ and $\tr A_jB_jC_j$,  are
bi-Lipschitz equivalent to their
comparison triangles with a constant  $L$ depending on $\lambda$ and $\delta$ only
(in notations of the lemma). This constant  can be chosen
as close to 1 as we wish, if  $\delta$ is small enough.
For a great $j$ both comparison triangles,
$\tr ABC$ и $\tr A_jD_jC_j$, are almost equal. So, all non-special triangles
of the surfaces $M_j$ are bi-Lipschitz equivalent to corresponding triangles of the
surface $M$.

Therefore, to finish the proof,  it is sufficient to verify
that (for great $j$) every special triangle
of the  surface $M_j$ is  bi-Lipschitz equivalent to the corresponding triangle
of the  surface $M$ or, equivalently, to its
comparison triangle.
\medskip

7. {\em Special triangles.}
Let $Q_0$ be the star of a fixed vertex
$B^0=E_k$ of the surfaces $M$,
$Q$ be the star of the corresponding
vertex $B=E_{kj}$ of the surfaces $M_j$. Recall that triangles of
 $Q_0$ are almost flat, so that they are
bi-Lipschitz equivalent to their
comparison triangles; the latter being  bi-Lipschitz equivalent to
comparison triangles for corresponding triangles of  $Q$.
This shows that it is sufficient to prove that (for sufficiently great $j$)
every star $Q$ is bi-Lipschitz equivalent to the star glued from
comparison triangles for the triangles of  $Q$.

8. {\em Plan of further proof.} We are going to apply to $Q$
 arguments from  \cite{BonkLang}. To do this  we attach a plane with a disk removed
to $Q$ and so we obtain a complete surface $P$ homeomorphic to the plane.
Recall that the key part
 of the proof in  \cite{BonkLang}
is, roughly speaking,  the following statement. If $P$ is a polyhedral surface
homeomorphic to the plane,
 $\omega^+(P)\leq2\pi-\epsilon<2\pi$, and
$\omega^-(P)\leq C<\infty$, then there is a set of flat
sectors with disjoint interiors on $P$;
 every point of nonzero curvature being a vertex for some  sectors.
We can decrease or increase (depending on the sign of curvature)
these  sectors so that
curvature at the sector vertices vanishes. Size of sector angles implies that this
process comes to a bi-Lipschitz map with a constant $L$ depending on $C$
and $\epsilon$ only. So we obtain a bi-Lipschitz map of $P$ to the plane $\R^2$.

Actually such a transformation requires three steps in \cite{BonkLang}.
First $P$ is divided onto two half planes by a special quasi-geodesic, and
the flat sectors are chosen separately in every half plane.
After that the vertices of positive curvature are removed. Finally
vertices of negative curvature are removed. See details in \cite{BonkLang}.

There is an obstacle for direct application of this construction in our case.
It is flat sectors containing rays that form small angles with the boundary
$\Gamma=\partial Q$ of
star $Q$.  Sectors on $P$ with vertices close to   $\Gamma$ can have such a property.
To avoid this difficulty, we choose  $j_0$ so great that almost all curvature
of $Q$ is concentrated in a very small neighborhood $V$ of the central point $B\in Q$
for $j>j_0$. After that, we replace a wide collar of $\partial Q$  by a flat
collar in $Q$. As a result, flat sectors come out to be almost orthogonal to
$\partial Q$ on the new deformed  surface. This simplifies further considerations.

9. {\em Elimination of curvature near   $\partial Q$.} Let $\{A_i\}$ be the set
of vertices of $\partial Q$, \;$|BA_i|=R_0$. Put
\begin{equation}
\label{dilatation}
\kappa=10\max\{\frac{2\pi}{2\pi-\epsilon},\,\frac{2\pi+C}{2\pi} \}.
\end{equation}
Consider the disk
${\bf D}(B,R)$ of radius $R$ such that
\begin{equation}
\label{short-begin}
10\kappa R<\delta R_0.
\end{equation}
After that, we choose  disk   ${\bf D}(B,r)$ (where $r\ll R$) and great number $j_0$
such that  $\mu^\pm (Q\setminus {\bf D}(B,r))<\delta$ for $j>j_0$ and, besides,
the conditions of Lemma
\ref{concentration} hold for $\Psi=C$, $Z=B$.

We are going to show that every star $Q$ is bi-Lipschitz equivalent to a region
$Q'$ (equipped with a polyhedral metric) which flat everywhere except a
$C^*r$-neighborhood of a point  $Z'$ located at a distance  $C^*R_0$
from the boundary of $Q'$;  $\mu^- (Q')< C+\delta$ and $\mu^+(Q')<2\pi-\frac12\epsilon$.

To simplify notations we omit indices and denote  by  $\tr BAC$ triangle
\linebreak
$\tr BA_iA_{i+1}$. Take points $A_1$, $C_1$ on the shortest curves $BA$, $BC$ at
distance $r/2$ from $B$ and connect these points by a shortest curve $A_1C_1$.

Let us strengthen our requirement about  $j_0$; namely, choose
$\rho>0$ so small and   $j_0$ so great that conditions of Lemma
 \ref{concentration}
 hold even if we replace $R$ and $r$ to $r$ and $\rho$, correspondingly.
In particular, we have  $\mu^\pm (Q\setminus {\bf D}(B,\rho))<\delta$.
It is not difficult to see that in this case   $A_1C_1$ is contained in
${\bf D}(B,r)$
and can not visit not only the disk ${\bf D}(B,\rho)$, but even
the disk ${\bf D}(B,r/4)$, and
angles  $\an BA_1C_1$ and $\an BC_1A_1$ are ``almost equal'' to angles
$\an A'_1$, $\an C'_1$ of comparison triangle $\tr B'A'_1C'_1$
(i.e., their differences are not greater than $\delta$).
In particular, these angles are less than
$\frac12 (\pi-\varphi_0)$. It is easy to see that the conditions
of Corollary  \ref{corr-4-angle} hold for the quadrangle
$AA_1C_1C$  (with an appropriate  $\phi$). Let us   apply the  corollary.
This allows us to replace each triangle $\tr BA_iA_{i+1}$ by a triangle
flat outside the disk ${\bf D}(B,r)$ and $L$-bi-Lipschitz
equivalent to $\tr BA_iA_{i+1}$.
Even if variation of  curvature of the new triangle is greater
than variation of  curvature of the old triangles
(at points  $A_1,\;C_1$), change of curvature
is not greater than $C^*\delta$. If we choose sufficiently small
$\delta$ and sufficiently
great $j_0$, we can take constant $L$ as close to 1 as we wish.

Let us save old notations $Q$, $B$, $A_1A_2\dots A_m$ for a new
star  arranged from the new  triangles and elements of the star.

Besides, we suppose $r$ to be so small in comparison  with $R$, that
$\tr ABC$ satisfies the conditions $A\in{\bf D}(B,\kappa r),\, |BC|\geq \kappa^{-1}R$.
\medskip

10. {\it Flat sectors.} We want to prove that a new star $Q$ is bi-Lipschitz
equivalent to a star obtained by gluing together comparison triangles for triangles
of the star $Q$.
To do this, we apply the construction from \cite{BonkLang}, described
above in short,  in item 8. This construction has to be applied twice:
first, to remove
positive curvature and, after that, to remove  negative one.
 This two steps are similar, so we will consider in details only the first one.

Let us supply $Q$ with a flat annulus to obtain an open
complete  surface $P$, flat everywhere except the disk ${\bf D}(B,r)\subset Q$.
This is possible. Indeed, denote by $\alpha^-_i$ and $\alpha^+_i$
adjacent to the base angles of triangle $\tr A_iBA_{i+1}$.
Consider a flat region bounded by two rays and interval of the length $|A_iA_{i+1}|$
under condition that angles between the interval and the rays from the region side
are equal to $\pi-\alpha^-_i,\; \pi-\alpha^+_i$, correspondingly.
Glue these flat regions together along rays and attach the  obtained region to $Q$.
For $j$ great enough, the  surface $P$ satisfies the conditions: its
positive curvature is less than $2\pi-\epsilon-\delta=2\pi-\epsilon'$ and negative
one is less than $C+\delta=C'$. Since our estimates are rough, we  preserve for
$\epsilon'$ and $C'$ previous notations $\epsilon$ and $C$.
\medskip

It follows from  \cite{BonkLang} that  there exists a finite set of flat sectors
with disjoint interiors on $P$ such that  all vertices of sectors are just
vertices of positive curvature  and the sum of angles for sectors with a common vertex $O$
equals
$$
\frac{2\pi-\omega^+(P)}{\omega^+(P)}\,\omega^+(O).
$$
To remove positive curvature at the point $O$,  we  stretch
all sectors with vertex $O$ by increasing their angles in
 $L_1=\frac{2\pi}{2\pi-\omega^+(P)}$ times. As a result, we obtain a polyhedron
$P_1$ of nonpositive curvature
bi-Lipschitz equivalent to $P$.

After this step, one finds an analogous system of flat sectors with vertices at
points of negative curvature and removes negative curvature  in the same way
by means of compressing  flat sectors of $P_1$. Finally we have a
bi-Lipschitz map $f\colon P_1\to\R^2$ with the Lipschitz constant
$\big(\frac{2\pi+C}{\epsilon})^\frac12$.

Following  \cite{BonkLang}, we use maps of the form  $(r, \phi)\to (r,a\phi)$
for  stretching and compressing sectors,  $(r,\phi)$ being polar coordinates with
origin  at the vertex of a sector.
We can assume that angles of the sectors are not big, in particular, that each sector
intersect only one special triangle base and the central point $B$ of the star $Q$
does not belong to the interior of a sector. To achieve this, it is enough to
divide sectors onto smaller ones.

Actually we will consider not all  surfaces $P$, but  only stars $Q$ of points $B$.
Such a star consists of isosceles triangles $\tr A_iB_jA_{i+i}$ and is bounded by
the geodesic broken $\Gamma=A_1A_2\dots A_m$. From description of the map $f$, it is
clear that $Q$ is bi-Lipschitz equivalent to a flat region --- its image $\tilde Q=f(Q)$.
The map $f$ transforms  bases  $\Gamma_i=A_iA_{i+1}$ of triangles
$\tr A_iB_jA_{i+i}$ to
curves  $\tilde\Gamma_i$ (not smooth in general). These curves consist of straight
segments (images of segments which do not belong to a flat sector) and
smooth curves (images of intersection of $\Gamma_i$ with a flat sector).
(It is not essential for us how
images of lateral  sides of triangles  $\tr A_iB_jA_{i+i}$
look like.)
 \medskip

11. {\em Flat region  $\tilde Q$.}  We are going to show that
the flat region  $\tilde Q$ is bi-Lipschitz
equivalent to a polygon
glued from comparison triangles for curved triangles of   $Q$.

Let us connect $\tilde B$ with points  $\tilde A_i$  by shortest curves in intrinsic
metric of  $\tilde Q$ (avoiding unnecessary intersections). So we divide  $\tilde Q$
onto ``curved triangles''  $\tilde T_i$ with curves $\tilde\Gamma_i$ as bases.
(It will be clear later that these shortest curves are almost orthogonal to
$\tilde Q$ and do not touch
one another.)

It is sufficient to verify that

(a) every curved triangle   $\tilde T_i$ is bi-Lipschitz equivalent to its
comparison triangle (i.e, a flat triangle $\tilde T'_i$
with side lengths equal to $|\tilde B\tilde A_i|$,
$|\tilde B\tilde A_{i+1}|$, and $s(\tilde\Gamma_i)$, correspondingly);

(b) the last flat triangle is  bi-Lipschitz equivalent to a comparison triangle
for  $\tr BA_iA_{i+1}$.

As the last triangle is almost equal (for great $j$) to the corresponding triangle
of the star $Q^0$, this ends the proof.

To prove (a) and (b), we need to estimate the  angle and
the distance distortions for map $f$.
To simplify exposition, we will consider only one step (removing positive curvature);
the estimates for the second step (removing negative curvature) are analogous.
\medskip

12. {\em Estimates.} To prove (a), we use Lemma  \ref{small_tr2}.
The following statements show that  $\tilde T_i$ satisfies the conditions of this lemma.
Also they help us to prove (b).

(i) For $j$ great enough, map $f$ slightly changes distances from $B$ to boundary $\Gamma$ of
 $Q$. More precisely, for every   $X\in \partial Q$ the inequality
\begin{equation}
\label{radial0}
||\tilde B\tilde X|-|BX||<C^*\delta |BX|
\end{equation}
holds.

(ii)  Let $X\in\tilde\Gamma_i$; then angles between radial shortest
curves $\tilde B\tilde X$ and
arcs of $\tilde\Gamma_i$ starting at $X$ are close
to  $\frac12 \pi$. In particular, turns of
 $\tilde\Gamma_i$ at its angular points are small. Words ``close'' and ``small''
 mean that difference between angles $\pi/2$ is not greater than $\frac{1}{10}$.

{\em Proof} (i).
  Let $X\in \partial Q$.
Prove that
$|BX|<(1+ \delta) |\tilde B\tilde X|$. The second required
inequality is proved  by analogy.

Consider a shortest curve  $\tilde\alpha$ connecting $\tilde B$ with
$\tilde X\in \partial \tilde Q$ and its
$f$-inverse image $\alpha$.  From \eqref{short-begin} it follows that the  initial
arc  $\tilde\alpha_0$ of $\tilde\alpha$, from $\tilde B$ to the boundary
of $f$-image of ${\bf D}(B, R)$,
is not longer than
$\kappa R<\frac{1}{10}\delta R_0\leq\frac{1}{5} \delta |BX|$.

If a segment of  the shortest curve $\alpha$ does not visit flat sectors,
map $f$ does not change  its length.
If a segment of  $\alpha$ is outside the disk  ${\bf D}(\tilde B, 2r_2)$ and
contained in a flat sector which was constricted,
it could become only shorter under $f$.

Now let  $\tilde\beta$ be an interval of   $\tilde\alpha$
containing in a flat sector $\tilde S$ such that $f$ got stretched $S$, and
$\beta$  $f$-inverse image in $S$ of the shortest path   $\tilde\beta$.
Let $O,\;\tilde O$ be vertices of sectors  $S,\,\tilde S$, correspondingly.

Denote by $Y,\;Z$  the initial and the end
points of segment
$\beta$, and by $\tilde Y,\;\tilde Z$  the initial and the end
points of  $\tilde\beta$. If $Z$ belongs to $\partial Q$, we replace
sectors $S$, $\tilde S$ by their subsectors  $ZOY, \,\tilde Z\tilde O\tilde Y$,
and preserve previous notations  $S,\;\tilde S$ for the new sectors.

We can suppose that $\tilde\beta$ does not intersect the initial segment
$\tilde\alpha_0$, so the distance between
$\tilde B$ and $\tilde Y$ is not less than $R$, and therefore (see
Lemma \ref{concentration}, item  (ii) and Remark \ref{remark})
$\an\tilde O\tilde Y\tilde B\leq 3\delta$.

Let us show that
\begin{equation}
\label{radial1}
s(\tilde\beta)\leq (1+C^*\delta)s(\beta),
\end{equation}
where, as usual,
$C^*$  means a constant depending on $C, \, \epsilon$ only.

Denote
$\an YOZ=\phi,\,\an \tilde Y\tilde O\tilde Z=\tilde\phi$,\,
$|OZ|=|\tilde O\tilde Z|=b$,\, $|YZ|=c$,\,
$|\tilde Y\tilde Z|=\tilde c$, $\pi -\an OYZ=\chi$.

 Place triangles  $\tr OYZ$ and $\tr\tilde S\tilde Y\tilde Z$
in $\R^2$ to one half-plane
with respect to  their common side $OY=\tilde O\tilde Y$.
Now it is clear that
$|\tilde c-c|\leq |\tilde ZZ|=2\sin\frac12(\tilde\phi-\phi)b\leq(\kappa-1)\delta c$,
as $\phi<\chi\leq\delta$. The last inequality follows from our choice
of disk ${\bf D}(B,r)$ in the beginning of item 9 and from
Lemma  \ref{concentration}.
So the estimate  \eqref{radial1}  is proved.

To obtain the second estimate it is enough to take the shortest curve $BX$
and its $f$-image in  capacity of $\alpha$ and $\tilde\alpha$ correspondingly.
\medskip

 {\em Proof} (ii).
We start with consideration of  $Q$ and, to be short,  denote $A_i=A$, $A_{i+1}=C$.
Side $AC$ is small in comparison with $|AB|=|CB|$ (see inequalities \eqref{size-angles}).
Therefore adjacent to base $A'C'$ angles of comparison triangle   $\tr A'B'C'$
are close to $\pi/2$ and angle  $\an A'B'C'$ is small.
The item  (i) of Lemma   \ref{concentration} says that angles
$\an BAC,\;\an BCA$ are close to $\pi/2$ either.

Consider a triangle  $\tr BAX$, where $X\in AC$. Let  $\tr B'A'X'$ be its
comparison triangle. Again from the item (i) of
Lemma \ref{concentration}, it follows that angles
$\an A, \;\an X$ are equal correspondingly to angles $\an A', \;\an X'$.
As $\an A'$ is almost equal $\pi/2$,
$\an A$ is close to  $\pi/2$ too. (``Close'' means that their difference
has the order of $0,01\epsilon+\delta$.)
Taking into account that angle $\an A'B'X'$ is small, from this it follows that
angle $\an A'X'B'$ is also close to $\pi/2$.
Now, again by item  (i) of Lemma  \ref{concentration}, angle  $\an AXB$ is close to
 $\pi/2$ too.
The same is true for angle $\an CXB$.

Let $S$ be a flat sector with a vertex $O$, sides of the sector intersect
$AC$ at points $X,\, Y$. Point $O$ is in the small neighborhood ${\bf D}(B,r)$ of $B$,
but not necessary in the triangle $\tr ABC$. It follows from  Lemma
\ref{concentration}, item (ii) that the angles
$\an OXB, \; \an OYB$ are small; therefore the angles
$\an OXY, \; \an OYX$  are close to  $\pi/2$
(by the same scale: their difference has the order $0,01\epsilon +\delta$).

Now pass to sector $\tilde S$, the image of flat sector $S$.
Radii of flat sector $S$
are almost orthogonal to
$\Gamma_i$.
A straightforward calculation shows that from this it follows that
radii of flat sector $\tilde S$
are almost orthogonal to $\tilde\Gamma_i$.  Distinction of the last angles
from $\pi/2$ depends on distinction between  angles $\Gamma_i$ and radii of flat
sectors $S$ from $\pi/2$ and on $\kappa$; i.e., finally on $C$ and $\epsilon$ only.

Vertices  $O$ of flat sectors $S$ are very close to $B$. Dilatation of $f$ is not
greater than $\kappa$, so $f$-images of vertices $O$ are close to $\tilde B$.
Hence,  $\an OXB$,  where $X\in \Gamma_i$, are close to zero, so angles between
segments $\tilde B\tilde X$ (they are shortest curves in $\tilde Q$)
and $\Gamma_i$    are close to  $\pi/2$.
(In particular, flat region $\tilde Q$ is a star region with respect to $\tilde B$.)
This proves item (ii).
\medskip

Estimate  (i)
implies that differences between  length of  the sides $\tilde A'\tilde B'$, \,
$\tilde C'\tilde B'$ of the
comparison triangle
$\tr \tilde A'\tilde B'\tilde C'$ and length of the sides
of the comparison triangle  $\tr A'B'C'$
are small. From   (ii) it follows that ratio of $|\tilde A'\tilde C'|$ to
 $|\tilde A'\tilde C'|$ is  bounded from below and above by numbers depending on
$C$ and $\epsilon$ only.
The choice of  $\varphi$ and item (ii) imply that
angle $\an \tilde A'\tilde B'\tilde C'$
is less than $\pi/2$. From this it becomes clear that flat triangles
$\tr \tilde A'\tilde B'\tilde C'$ and $\tr A' B' C'$ are $L$-bi-Lipschitz equivalent,
where $L$ depends on $C$ and $\epsilon$ only; for example, see Corollary  1
in the paper  \cite{BeBu}.
Finally, each triangle  $\tr \tilde A\tilde B\tilde C$
is bi-Lipschitz equivalent to the corresponding triangle of the star
$Q^0$, and our theorem is proved completely.

\end{document}